\documentclass{amsart}
\usepackage{import}
\usepackage[latin9]{inputenc}
\usepackage{microtype}
\usepackage{fancyhdr}
\usepackage{cite} 
\usepackage{amsfonts}
\usepackage{amsmath}
\usepackage{amsthm}
\usepackage{amssymb}
\usepackage[all]{xy}
\usepackage{mathrsfs}
\usepackage{bm} 
\usepackage{stmaryrd}
\usepackage{color}   
\usepackage[hidelinks]{hyperref}
\hypersetup{
    colorlinks=false, 
    linktoc=all,     
    linkcolor=black,  
}
\usepackage{graphicx}
\usepackage{wrapfig}
\usepackage{float}
\usepackage{caption}
\usepackage{subcaption}
\usepackage{enumitem}
\usepackage{empheq} 
\usepackage{mathtools} 
\usepackage{adjustbox}
\usepackage{tikz,tikz-cd}
\usetikzlibrary{positioning,arrows}
\setlist{nosep} 
\usepackage{dutchcal} 

\usepackage{makecell} 

\usepackage[foot]{amsaddr} 

\usepackage{imakeidx} 
\makeindex 

\usepackage[splitrule]{footmisc} 

\usepackage[top=2.5cm,bottom=2cm, margin=2.5cm]{geometry}
\usepackage[font=small, width=.75\textwidth]{caption}
\usepackage{mwe}

\theoremstyle{definition}
\newtheorem{theorem}{Theorem}

\theoremstyle{definition}

\theoremstyle{definition}

\theoremstyle{definition}
\newtheorem{proposition}[theorem]{Proposition}

\theoremstyle{definition}

\theoremstyle{definition}

\theoremstyle{definition}

\theoremstyle{definition}

\newtheorem*{innerthm}{\innerthmname}
\newcommand{\innerthmname}{}
\newenvironment{statement}[1]
 {\renewcommand{\innerthmname}{#1}\innerthm}
 {\endinnerthm}
\theoremstyle{definition}

\usepackage{titlesec}
\titleformat{\section}
  {\normalfont\fontsize{14}{15}\bfseries}{\thesection}{1em}{}

\titleformat{\subsection}
  {\normalfont\fontsize{12}{15}\bfseries}{\thesubsection}{1em}{}

\usepackage{xcolor}

\newcommand*{\colorboxed}{}
\def\colorboxed#1#{%
  \colorboxedAux{#1}%
}
\newcommand*{\colorboxedAux}[3]{%
  \begingroup
    \colorlet{cb@saved}{.}%
    \color#1{#2}%
    \boxed{%
      \color{cb@saved}%
      #3%
    }%
  \endgroup
}

\title{\vspace{-3cm}Obstructing Anosov Flows on Cusped 3-Manifolds}
\author{Misha Schmalian}
\address{\hskip -0.35cmMathematical Institute, University of Oxford, OX2 6GG, UK, https://sites.google.com/view/mishaschmalian}
\begin{document}
\vspace*{3cm}
\maketitle
\date{}
\vspace{-1cm}
\begin{abstract}
Using results relating taut foliations and pseudo-Anosov flows, we find cusped hyperbolic 3-manifolds which are not the non-singular part of a pseudo-Anosov flow. In particular, we find the first examples of cusped hyperbolic 3-manifolds not admitting veering triangulations, confirming a conjecture of S.  Schleimer.
\end{abstract}

\section{Introduction}

\noindent
Consider a closed orientable 3-manifold $M$. Anosov flows, defined below, present one of the most well-behaved and well-studied classes of dynamical systems on $M$. For example, the geodesic flow of hyperbolic surfaces and the suspension flow of hyperbolic toral homeomorphisms are both Anosov. The question of which manifolds admit Anosov flows was raised early on and obstructions due to group growth, foliations, and contact topology have been found, see \cite[Appendix]{Margulis_Anosov_Growth}, \cite{RobertsShareshianStein}, \cite{Etnyre_Contact_Anosov}. In this article, we aim to obstruct Anosov flows on cusped hyperbolic manifolds. However, we have to consider Anosov flows with nice boundary behaviour. \\

\noindent
A pseudo-Anosov flow on $M$, also defined below, is an Anosov flow away from a finite link $\Gamma\subset M$ of singular orbits. We can think of a pseudo-Anosov flow as an Anosov flow on $M-\Gamma$ with particularly good behaviour at the cusps. Call a cusped 3-manifold $M$ {\it persistently foliar} if there is a choice of slope $l_i$ in each boundary component, s.t. a Dehn filling $M(s_1, \ldots, s_k)$ admits a coorientable taut foliation whenever $s_i\neq l_i$ for all $i=1, \ldots, k$. The following obstruction to pseudo-Anosov flows follows from techniques well-known to experts: 

\begin{proposition} \label{anosov_foliar} Given a closed 3-manifold $M$ admitting a pseudo-Anosov flow with singular orbits $\Gamma$, some double cover of $M-\Gamma$ is persistently foliar. 
\end{proposition}

\noindent
A combinatorial analog of pseudo-Anosov flows are veering triangulations, which again we define below. These were introduced in \cite{Agol_veering} and can be used as a complete conjugacy invariant for pseudo-Anosov maps of surfaces. A veering triangulation on a cusped 3-manifold $M$ gives a pseudo-Anosov flow on a Dehn-filling $N=M(s_1, \ldots, s_k)$, provided each $s_i$ avoids a specific slope on each component of $\partial M$. This has singular orbits $\Gamma$, with $N - \Gamma \cong M$, and without, so called, perfect fits \cite{ChiCheukVeering, Veering_DynamicsPairs}. Conversely, such pseudo-Anosov flows on $N$ give rise to veering triangulations on $N-\Gamma$ \cite{AgolGueritaudSummary, Veering-LoomSpaces} and moreover, by work in preparation \cite{Veering_Pseudo_Anosov_Back}, these maps give a 1-1 correspondence. Veering triangulations also relay geometric information. They admit strict angle-structures \cite{VeeringAngleStructure, Angle_Structures_Futer} therefore forcing the underlying manifold to be hyperbolic \cite[Corollary 4.6.]{Angle_Structures}. It is thus instructive to ask whether every cusped hyperbolic manifold admits a veering triangulation, and in fact S. Schleimer has conjectured this is not the case \cite{SaulConjectureVeering}. The correspondence between veering triangulations and pseudo-Anosov flows without perfect fits gives the following: 

\begin{proposition} \label{veering_foliar}
Any cusped manifold admitting a veering triangulation has a persistently foliar double cover. 
\end{proposition}

\noindent
There is an obstruction theory for coorientable taut foliations using the Heegaard-Floer homology of the manifold \cite{OsvathSzaboLSpace}. This was used to create a large data set of 3-manifolds not admitting taut foliations \cite{DunfieldData}. Using the computer programs SnapPy \cite{SnapPy} and this data set we can find cusped hyperbolic manifolds without persistently foliar double covers. This gives the following:

\begin{theorem} \label{no_veering}
There exists cusped hyperbolic 3-manifolds which do not admit veering triangulations. Moreover, these cusped 3-manifolds are not the non-singular part of a pseudo-Anosov flow.
\end{theorem}

\noindent
Examples of manifolds not admitting veering triangulations include the Whitehead link complement and m006 of the cusped census. In fact, among the first 100 manifolds in the orientable cusped census we show 19 have no veering triangulations, compared to the 20 appearing in the veering census \cite{VeerinCensus}. \textcolor{black}{We should mention that an alternate proof that the Whitehead link complement admits no veering triangulations is in preparation \cite{UniqueVeering}.} To be clear, whenever we show a cusped manifold has no veering triangulations, we also show this manifold is not the non-singular part of a pseudo-Anosov flow. In particular, we do not make use of the lack of perfect fits in the pseudo-Anosov flow of a veering triangulation. 

\section{Acknowledgements} 
\noindent
The author would like to thank Marc Lackenby for all his time, advice and support. The author would also like to thank Henry Segerman and Saul Schleimer for helpful discussion and feedback on this article.
\section{Definitions}
\noindent
Let us deliver the promised definitions. 

\noindent
An {\it Anosov flow} $\Phi_t$ on an orientable closed 3-manifold $M$ is a flow preserving the splitting $TM=E^s\oplus E^u\oplus X$, s.t. the flow is tangent to $X$ and there are $C\geq 1, \lambda>0$ with: 
$$||D\Phi_t(v)||\leq C e^{-\lambda t} ||v||, ||D\Phi_{-t}(w)||\leq C e^{-\lambda t} ||w||, \text{ for any }v\in E^s, w\in E^u, t\geq 0.$$
The plane bundles $TX\oplus E^s$ and $TX\oplus E^u$ are integrable and give the {\it (weak) stable foliation} $\mathcal{F}^s$ and the {\it (weak) unstable foliation} $\mathcal{F}^u$ of the flow respectively. A foliation is {\it taut} if every leaf is intersected by a loop transverse to the foliation. The foliations $\mathcal{F}^s, \mathcal{F}^u$ of an Anosov flow are known to be taut.\footnote{This is because leaves of the stable foliations are contracted/dilated, hence cannot be compact, and in particular do not bound a dead-end component.}. \\

\noindent
A flow $\Phi_t$ on $M$ is {\it pseudo-Anosov} if it is Anosov away from a finite set of closed orbits, called the {\it singular orbits}, and on the singular orbits $\Phi_t$ is locally modeled on the suspension flow of a pseudo-Anosov diffeomorphism. A closed manifold $M$ admitting a pseudo-Anosov flow with singular orbits $\Gamma$ can be thought of as an Anosov flow on $M-\Gamma$ with particularly well-behaved behaviour at the boundary. \\

\noindent
{\it Veering triangulations} are ideal triangulation $\mathcal{T}$ of cusped 3-manifolds with each 2-simplex of $\mathcal{T}$ transversely oriented and each 1-simplex of $\mathcal{T}$ coloured red/blue, s.t. the following conditions are satisfied: 
\begin{itemize}
\item[-] At each 1-simplex of $\mathcal{T}$ at most two pairs of consecutive 2-simplices are not compatibly oriented; 
\item[-] Within each 3-simplex of $\mathcal{T}$ two 2-simplices are oriented inwards, and two outwards; 
\item[-] At each 3-simplex, the edges between compatibly oriented 2-simplices are coloured as shown below. 
\end{itemize}
\noindent
Note: The first two conditions intuitively say the 3-simplices are ``flat''. Moreover, in the figure below the two remaining uncoloured black edges can be coloured red or blue depending on the rest of the triangulation.
\vskip 0.2cm
{\centering
\includegraphics[width=3.8cm]{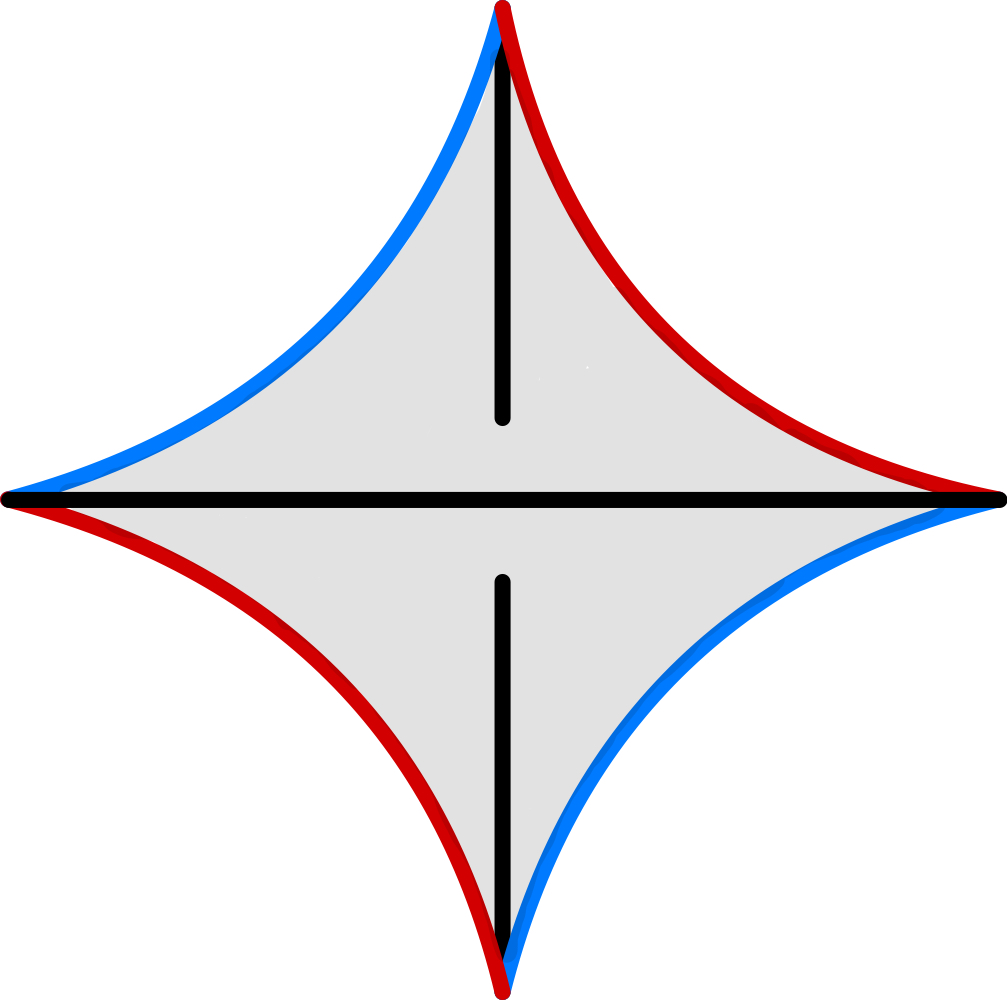}
\captionof{figure}{3-simplex of a veering triangulation, with transverse orientations of 2-simplices pointing out of the page.}
\label{veering_oriented}
}
\vskip 0.2cm
\noindent
Recall that we say a cusped 3-manifold $M$ is {\it persistently foliar} if there is a choice of slope $l_i$ in each boundary component, s.t. a Dehn filling $M(s_1, \ldots, s_k)$ admits a coorientable taut foliation whenever $s_i\neq l_i$ for all $i$.

\section{Proof of Theorem}

\noindent
The following result is well-known to experts:

\begin{theorem}\label{Main_theorem} Consider a compact 3-manifold $M$ with toroidal incompressible boundary and taut coorientable foliation $\mathcal{F}$. If $\mathcal{F}$ intersects each boundary torus in a codimension-1 foliation consisting only of closed curves and a non-zero number of Reeb annuli, then $M$ is persistently foliar. 
\end{theorem}

\noindent {\it Proof sketch}: \cite[Example 4.22]{Calegari_foliation_book} \cite[Operation 2.4.1]{Gabai_Suspension_S1} At a boundary torus $T\subset \partial M$ glue in a solid torus $S$ along any slope which is not given by the closed curves of $\mathcal{F}\cap T$. Within $S$ cap of the Reeb annuli with annular leaves and extend the foliation to the solid tori they bound. Blow up the leaves $\lambda_i$ meeting $T$ in closed curves replacing them with $I$-bundles $\lambda_i\times I$. The yet unfoliated region of $S$ is an ideal polygonal bundle with the polygon fibre having an even number of vertices due to coorientability of $\mathcal{F}$. We foliate this region with a family of ``stacked''  monkey saddles. Finally we need to extend the foliation of $(\partial \lambda_i)\times I$ given by the monkey saddles to all of $\lambda_i\times I$. This can be done by finding a holonomy homomorphism $\rho:\pi_1(\lambda_i)\to \text{Homeo}^+(I)$ which restricts to the desired holonomy on $\partial \lambda_i$. The existence of such a holonomy is detailed in \cite{Calegari_foliation_book}. Denote the new foliation by $\mathcal{F}'$. \\
Every leaf of $\mathcal{F}'$ was already a leaf of $\mathcal{F}$ or enters the blow up $\lambda_i\times I$ of a leaf. Thus a loop $\gamma\subset M$ transversely intersecting every leaf of $\mathcal{F}$ will continue to transversely intersect every leaf of $\mathcal{F'}$ and $\mathcal{F'}$ is taut. \\
Since $\rho$ takes values in $\text{Homeo}^+(I)$, $\mathcal{F'}$ is orientable away from the stacked monkey saddles. In $\mathcal{F}\cap T$ adjacent Reeb components are oriented in opposite direction. This is consistent with the boundary of a monkey saddle alternating between ``upwards bending'' and ``downwards bending'' leaves. Thus $\mathcal{F'}$ is coorientable and taut and $M$ is persistently foliar. ``$\square$''\\

{\centering
\includegraphics[width=13.2cm]{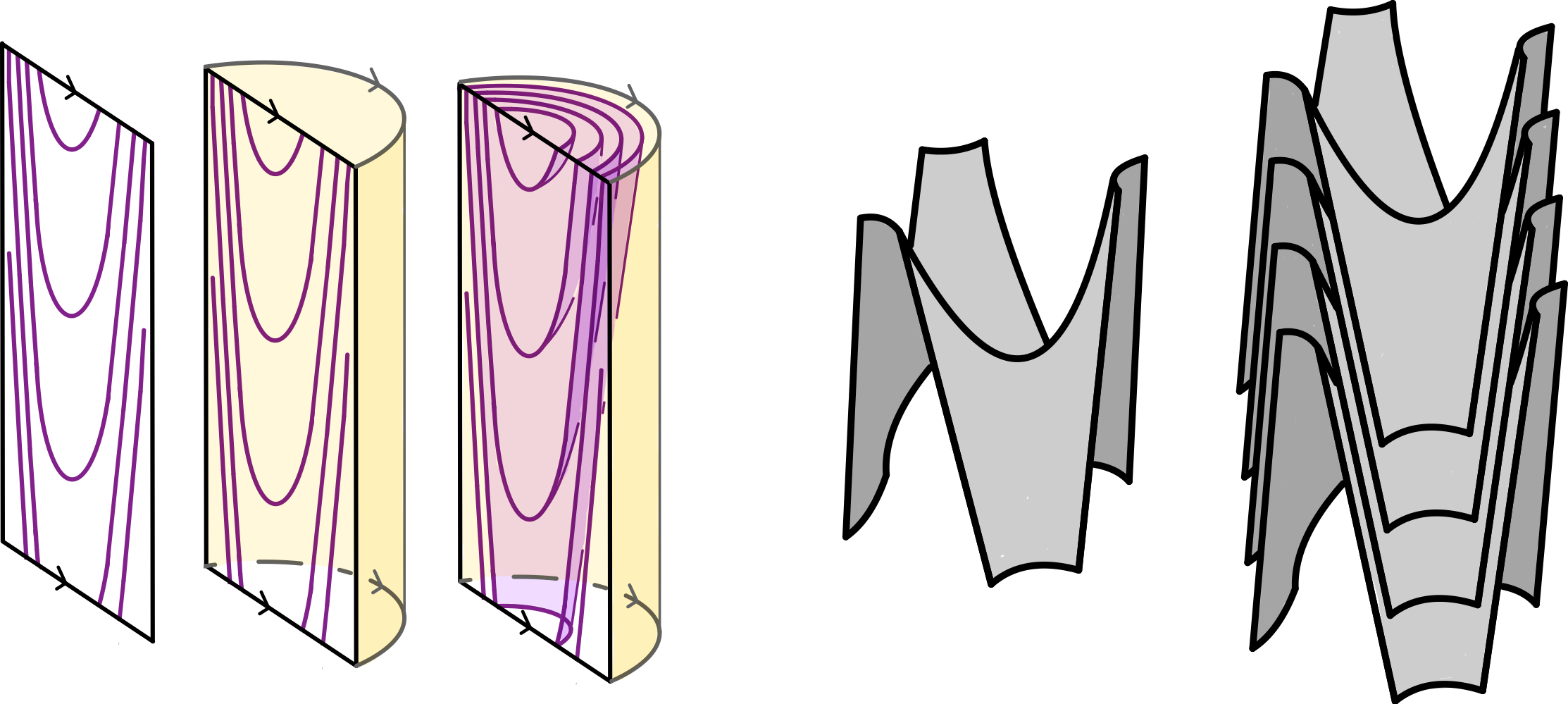}
\captionof{figure}{Left: Annuli Reeb foliation in $\mathcal{F}\cap T$ in purple, capped off by an annular leaf in yellow, and an extension of $\mathcal{F}$ to the solid tori in purple again. Right: A monkey saddle next to a stack of monkey saddles giving a foliation of an ideal polygon bundle.}
\label{mapping_torus}}
\vskip 0.3cm

\noindent
Let us now prove the previously stated result:

\vskip 0.2cm
{\centering
\includegraphics[width=4.7cm]{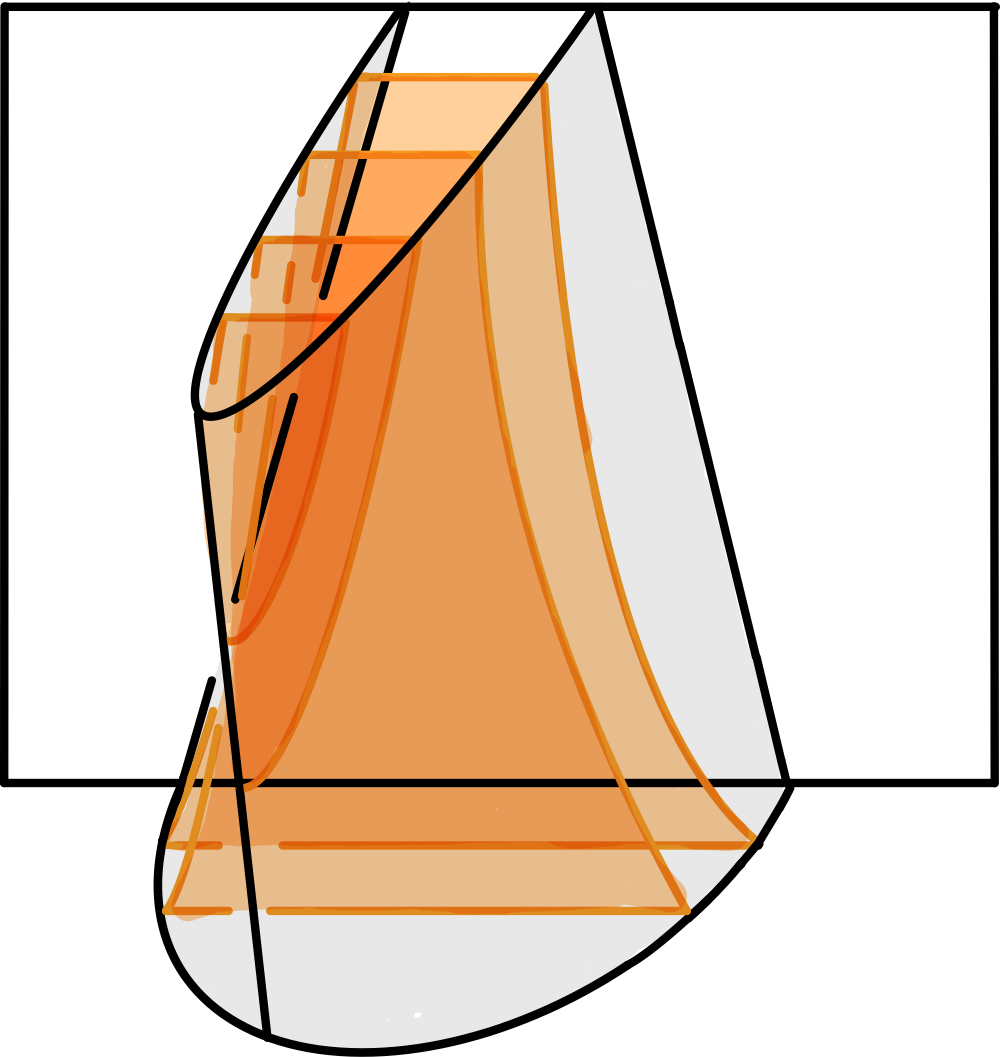}
\captionof{figure}{Suspension of half disc with boundary shaded in grey and foliation by leaves $\mathbb R\times \{p\}\times \mathbb R$ in orange.}
\label{Anosov_folio}}
\vskip 0.3cm

\begin{statement}{Proposition 1}
Given a closed 3-manifold $M$ admitting a pseudo-Anosov flow with singular orbits $\Gamma$, some double cover of $M-\Gamma$ is persistently foliar. 
\end{statement}

\noindent
{\it Proof}: Let $\mathcal{F}^s$ denote the stable foliation of the Anosov flow on $M-\Gamma$. It is sufficient to find a tubular neighbourhood $N$ of $\Gamma$ with $\mathcal{F}^s\cap \partial N$ being foliated by Reeb annuli and closed curves. In this case, the double cover of $M-N$ for which the lift of $\mathcal{F}^s$ is coorientable satisfies the criteria of Theorem \ref{Main_theorem} and therefore the double cover is persistently foliar. \\
\noindent
We may take $N$, s.t. $\partial N$ transversely intersects the singular leaf $\lambda_s$ of $\mathcal{F}^s$ in closed curves and s.t. components of $N-\lambda_s$ are modelled on the suspension of an Anosov map $\phi: \mathbb R\times \mathbb R_{>0}, \to \mathbb R\times \mathbb R_{>0}$ given by $\begin{pmatrix} a & 0\\ 0& 1/a \end{pmatrix}, a>1.$ We may further shrink $N$, s.t. each component of $N-\lambda_s$ is the suspension of a half disc: $\left\{\left(\begin{pmatrix} a^t & 0 \\ 0 & a^{-t}\end{pmatrix}\cdot z, t\right): z\in \mathbb R\times \mathbb R_{>0}, |z|<1, t\in [0,1]\right\}\subset\left( \mathbb R\times \mathbb R_{>0}\times [0,1]\right)/\phi$. Under this identification $\mathcal{F}^s$ consists of vertical planes and foliates components of $\partial N-\lambda_s$ as Reeb annuli, see figure \ref{Anosov_folio}. $\square$

\begin{statement}{Proposition 2}
Any cusped manifold admitting a veering triangulation has a persistently foliar double cover. 
\end{statement}

\noindent
This follows since any manifold with veering triangulation is the non-singular part of a pseudo-Anosov flow.

\noindent
We call a manifold {\it taut} if it admits a coorientable taut foliation. There is a large data set of 3-manifolds known to be, or not to be, taut \cite[floer/final\underline{~}data/QHSpheres]{DunfieldDataSpecific}. We thus proceed as follows: start with a cusped hyperbolic manifold $M$; use SnapPy to identify the double covers $N_1, \ldots, N_k$ of $M$; use the data from \cite{DunfieldDataSpecific} to find sufficiently many non-taut fillings of each of the $N_i$ to show that $N_i$ are not persistently foliar; apply Propositions \ref{anosov_foliar} and \ref{veering_foliar}.\\

\noindent
We fix a basis convention to describe slopes for Dehn-filling. For once-cusped manifolds we adopt the convention chosen in \cite{DunfieldDataSpecific}. For multi-cusped manifolds we use the basis selected by SnapPy. The exact convention is not critical if we are consistent when working with a fixed manifold.\\

\noindent
For example, the once-cusped manifold m006 has a single double cover s649 and the Dehn-fillings s649(-1/3) and s649(-1/4) do not admit taut foliations. Therefore m006 is not the non-singular part of a pseudo-Anosov flow, and in particular does not admit a veering triangulation. This proves Theorem \ref{no_veering}.\\

\noindent
For another example, consider the Whitehead link L5a1. Its complement has two double covers, namely the complement of the three component link L8n5, shown in figure \ref{Whitehead_link}, and the twice cusped t12048. 

\vskip 0.15cm
{\centering
\includegraphics[width=8cm]{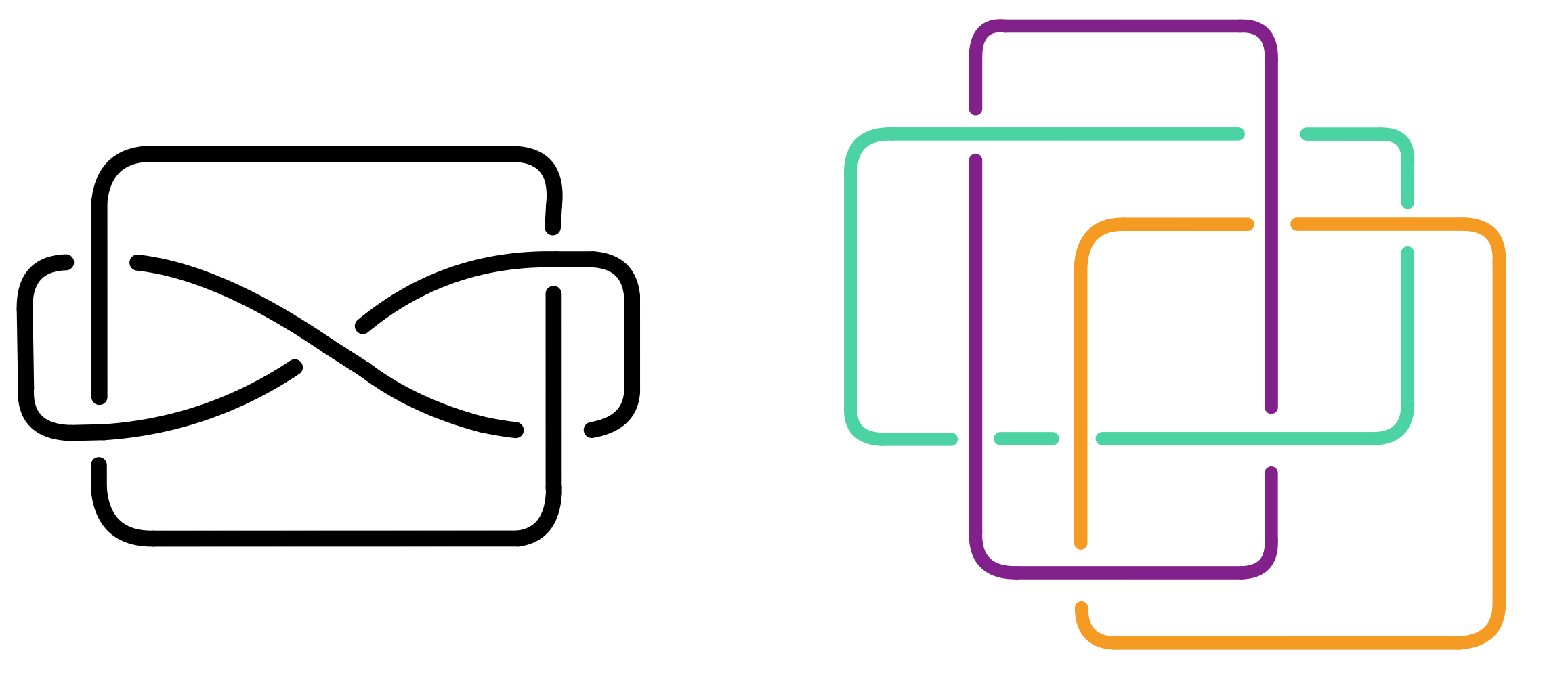}
\captionof{figure}{Whitehead link L5a1 (left) and the Link L8n5 (right) with cusp 0 in yellow, cusp 1 in green, and cusp 2 in purple.}
\label{Whitehead_link}
}
\vskip 0.3cm

\noindent
Observe that $\text{L8n5}(*, \infty, *)\cong \text{L8n5}(*, *, \infty)$ is the complement of the two-component unlink, and in particular does not admit a taut foliation \cite{Novikov_Theorem}. Thus L8n5 can only be persistently foliar with exceptional slopes $l_1, l_2=\infty$. SnapPy identifies that $\text{L8n5}(* ,1/2, 2)\cong \text{m149}$ and the Dehn-fillings m149(-1/3), m149(-1/4) are non-taut by \cite{DunfieldDataSpecific}. Thus L8n5 is not persistently foliar. SnapPy identifies that $\text{t12048}(2, *)\cong \text{s778}$ and $\text{t12048}(3, *)\cong \text{t07936}$. The Dehn-fillings s778(-3/1), s778(-4/1), t07936(-4/1) and t07936(-5/1) are non-taut by \cite{DunfieldDataSpecific}. This shows that t12048 is not persistently foliar and hence that the Whitehead link is not the non-singular part of a pseudo-Anosov flow nor admits a veering triangulation.\\

\noindent
We mention a criteria for a once-cusped 3-manifold to not be persistently foliar that is not directly used here but is still of interest. An {\it L-space} is a closed rational homology 3-sphere $M$ with Heegaard-Floer homology satisfying $\dim\widehat{HF}(M)=|H_1(M; \mathbb Z)|$. A knot $K$ in an L-space $Y$ {\it Floer simple} if $\widehat{HFK}(K)\cong \mathbb Z^{|H_1(Y-K)|}$ and a once-cusped manifold is {\it Floer simple} if it is the complement of a Floer simple knot. Once-cusped Floer simple 3-manifolds admit multiple L-space fillings \cite[Proposition 1.3]{Rasmussen_Rasmussen_Floer_Simple}. Thus if every double cover of a cusped 3-manifolds $M$ is once-cusped and Floer simple, then $M$ is not the non-singular part of a pseudo-Anosov flow, nor admits veering triangulations. One could hope for this sort of approach to obstruct pseudo-Anosov flows. \\

\noindent
\textcolor{black}{A veering triangulation is {\it edge-orientable} if one may orient the edges in a manner satisfying certain conditions, see \cite[Definition 8.1.6]{Parlak_Thesis}. Edge-orientable veering triangulations are exactly those that give rise to Anosov flows with coorientable stable and unstable foliations. Similarly to the work above, we obtain a simplified obstruction for edge-orientable veering triangulations.}
\begin{statement}{Proposition 3}
Any cusped manifold admitting an edge-orientable veering triangulation is persistently foliar. 
\end{statement}

\section{More Examples}
\noindent
We repeat the technique used on m003 and the Whitehead link complement for the first 100 manifolds of the orientable cusped census and tabulate results below. 

\begin{center}
\begin{tabular}{ | c | }
\hline
Appears in the veering census:\\
\hline 
\makecell{m003, m004, m009, m010, m016, m022, m023, m036, m038, m039, \\m040, m052, 
m083, m115, m119, m120, m125, m135, m136, m140}\\
\hline\hline
Has no veering triangulations: \\
\hline
\makecell{m006, m007, m011, m029, m030, m037, m047, m049, m060, m064,\\ m081, m082, m095, m116, m117, m129, m130, m142, m143}\\
\hline
\end{tabular}
\end{center}

\noindent
For the census manifolds above we tabulate the census name of their double covers and a sufficient number non-taut Dehn-fillings of the double covers. As before, double covers are identified via SnapPy and non-taut fillings are found using \cite{DunfieldDataSpecific}. 

\begin{center} 
\begin{tabular}{| c || c | c | c | c | c |}
\hline
Manifold & m006 & m007 & m011 & m029 & m030\\
\hline 
Double Covers & s649 & s649 & s874 & t07681  & t07681 \\
\hline
Initial Fillings & & & & & \\
\hline
Non-taut Fillings & \makecell{s649(-1/3), \\ s649(-1/4).} & \makecell{s649(-1/3). \\ s649(-1/4).} & \makecell{s874(-1/3), \\ s874(-1/4).} & \makecell{ t07681(1/5), \\  t07681(2/5).}  & \makecell{ t07681(1/5), \\  t07681(2/5).}\\
\hline
\end{tabular}
\end{center}
\begin{center}
\hskip 0.45cm
\begin{tabular}{ || c | c | c | c | c | } 
\hline
  m037 & m047 & m049 & m060 &  m064  \\
\hline 
  t07933, t07939, v3222 & v3387 & v3431 & t09618 &  t09795 \\
\hline
  \makecell{v3222(1$,*$) $\cong$ m035,\\ v3222(2$,*)\cong$ m307.} & & &  &  \\
\hline
  \makecell{t07933(-1/3), t07933(-1/4),\\ t07939(1/5), t07939(2/5),\\ m035(-2/3), m035(-4/3), \\m307(-1/3), m307(-1/4).
} & \makecell{v3387(-1/4), \\v3387(-1/5).} & \makecell{v3431(-1/4), \\v3431(-1/5).} & \makecell{t09618(-1/4), \\t09618(-1/5).} &  \makecell{t09795(1/4), \\t09795(1/5).}  \\
\hline
\end{tabular}
\end{center}

\begin{center}
\hskip -1.68cm
\begin{tabular}{  || c | c | c | c | c | } 
\hline
m081 & m082 & m095 & m116 &  m117 \\
\hline 
 t10615 & t10708 & t10831 & t11579 &  t11693 \\
\hline
 & & &  &  \\
\hline
\makecell{t10615(-4), \\ t10615(-5).}  & \makecell{t10708(-1/4), \\t10708(-1/5).} & \makecell{ t10831(-1/5), \\t10831(-5).} & \makecell{t11579(-1/3), \\t11579(-1/4).} &  \makecell{t11693(-1/3), \\t11693(-1/4).}  \\
\hline
\end{tabular}
\end{center}

\begin{center}
\hskip -0.94cm
\begin{tabular}{|| c | c | c | c | c |} 
\hline
m129 & m130 \\
\hline 
t12066, t12048 & t12038, t12048 \\
\hline
  \makecell{t12066$(2, 2, *)\cong$ m149,  t12066$(3, 3, *)\cong$ s673, \\ t12066$(2, 3, *)\cong$m307,  t12066$(3, 2, *)\cong$ m288,\\  t12048$(2, *)\cong$ s778, t12048$(3/2, *)\cong$ t08875.}  &\makecell{t12048$(2,*)\cong$ s778,\\ t12048$(3/2, *)\cong$ t08875.} \\
\hline
\makecell{m149(-1/3), m149(-1/4), s673(-1/3), s673(-1/4), \\m307(-1/3), m307(-1/4), m288(1/3), m288(1/4), \\s778(-3), s778(-4), t08875(-1/4), t08875(-2/3).}&\makecell{t12038(-1/3), t12038(-1/4), \\s778(-3), s778(-4), \\t08875(-1/4), t08875(-2/3).} 
\\
\hline
\end{tabular}
\end{center}

\begin{center}
\hskip-8.26cm
\begin{tabular}{|| c | c | } 
\hline
 m142 & m143 \\
\hline 
 t12310  & o9\underline{~}33110\\
\hline
&  \\
\hline
 \makecell{t12310(-3,1), \\t12310(-4).} & \makecell{o9\underline{~}33110(-4), \\o9\underline{~}33110(-4/3).}
\\
\hline
\end{tabular}\\
\end{center}

\vskip 0.4cm
\noindent
\textcolor{black}{We may use Proposition 3 and \cite{DunfieldDataSpecific} to obstruct edge-orientable veering triangulations. For example, we see that m003 admits a veering triangulation, but is not  persistently foliar and thus does not admit an edge-orientable veering triangulation. In fact, among the 1426 manifolds of the veering census that are in the cusped census, 503 do not admit edge-orientable veering triangulations. Moreover, \cite{DunfieldDataSpecific} gives 50598 hyperbolic once-cusped 3-manifolds, which are Floer simple and hence do not admit edge-orientable veering triangulations.}

\bibliographystyle{alpha}
\bibliography{citation}

\end{document}